\documentclass[final]{jotart}
 
\usepackage{amsfonts}
\usepackage{amssymb,amsmath}
\usepackage{amsthm,amscd}
\usepackage{upref}
\usepackage{enumerate}
\usepackage{bbm}
 
\theoremstyle{proclaim}
\newtheorem{theorem}{Theorem}[section]
\newtheorem{lemma}[theorem]{Lemma}

\newtheorem{proposition}[theorem]{Proposition}

\theoremstyle{statement}
\newtheorem{remark}[theorem]{Remark}

\newtheorem{example}[theorem]{Example}

\newcommand{\1}{\mathbbm{1}}

 \newcommand{\supp}{\operatorname{supp}}
  \newcommand{\e}{\eqref}

\newcommand{\q}{\quad}

\newcommand{\wh}{\widehat}
\newcommand{\la}{\langle}
\newcommand{\ra}{\rangle}

\newcommand{\ov}{\overline}
 \renewcommand{\d}{\delta}
 
 \newcommand{\cl}{\operatorname{clos}}

\newenvironment{pf}{\begin{proof}}{\end{proof}}

\def\qqq{\mathrel{\subset\mkern-15mu\lower.38ex\hbox{${\scriptscriptstyle\rightarrow}$}}}

\let\cal\mathcal

\let\Bbb\mathbb

\numberwithin{equation}{section}

\begin{document}

 \issueinfo{vv}{n}{yyyy} 
\commby{A. B\"ottcher}
\pagespan{000}{000}
\date{Month dd, yyyy}
\revision{Month dd, yyyy}

\title {On semibounded Toeplitz  operators}

\author{ D. R. Yafaev}
\address{ IRMAR, Universit\'{e} de Rennes I\\ Campus de
  Beaulieu, 35042 Rennes Cedex, FRANCE}
\email{yafaev@univ-rennes1.fr}

\begin{abstract} 
 We show that a semibounded Toeplitz quadratic form   is closable in the space $\ell^2({\Bbb Z}_{+})$ if and only if
 its entries  are Fourier coefficients of an absolutely continuous measure.  We also describe the domain of the corresponding closed form. This allows us to define semibounded Toeplitz operators under minimal assumptions on their matrix elements. 
     \end{abstract}
     
      \begin{subjclass}{Primary 47A05, 47A07; Secondary 47B25, 47B35}
      \end{subjclass}


\maketitle


\section{INTRODUCTION. MAIN RESULTS }  

{\bf 1.1.}
Toeplitz operators   $T$ can formally be defined in the space $\ell^2({\Bbb Z}_{+})$ of sequences $g=(g_{0}, g_{1}, \ldots)$
 by the formula
  \begin{equation}
(T g)_{n}=\sum_{m=0}^\infty  t_{n-m} g_m, \q n=0,1,\ldots.
 \label{eq:HD}\end{equation}
 Thus the matrix elements of a  Toeplitz operator depend  on the difference of the indices only. So it is natural to expect that properties of  Toeplitz operators are close to those of discrete convolution operators acting in the space $\ell^2({\Bbb Z})$.
  
 The precise definition of the operator   $T$ requires some accuracy.  Let  $\cal D \subset \ell^2({\Bbb Z}_{+})  $ be the dense set   of sequences $g=\{g_{n}\}_{n\in {\Bbb Z}_{+}} $ with only a finite number of non-zero components. If the sequence 
 $t=\{t_{n}\}_{n\in {\Bbb Z}}  \in \ell^2({\Bbb Z})  $, then for $g\in \cal D $ sequence \e{eq:HD}   belongs to $ \ell^2({\Bbb Z}_{+})  $. In this case the operator $T$ is defined on $\cal D $, and it is symmetric  if $t_{n}= \ov{t_{-n}}$. Without any a priori assumptions on $t_{n}$, 
  only the quadratic form
   \begin{equation}
t[g,g] =\sum_{n,m\geq 0} t_{n-m} g_{m}\bar{g}_{n}   
 \label{eq:QFq}\end{equation}
 is well defined for $g\in \cal D$.
 
 The theory of Toeplitz  operators is a very well developed subject. 
 We refer to the books \cite{Bo} and \cite{GGK} (Chapter~XII),  \cite{NK} (Chapters~B.4 and B.6),   \cite{Pe} (Chapter~3)  for basic information on this theory.

   Let us state a necessary and sufficient condition  for a Toeplitz operator $T$   to be bounded. Below $d{\bf m}$ is always the normalized Lebesgue measure on  the unit circle $\Bbb T$.

 \begin{theorem}[Toeplitz]\label{B-H} 
  A Toeplitz operator $T$ $($defined, possibly, via its quadratic form \e{eq:QFq}$)$ is bounded if and only if the $t_{n}$ are the Fourier coefficients of some bounded function on  $\Bbb T$:
    \begin{equation}
  t_{n} =\int_{\Bbb T} z^{-n} w (z) d{\bf m} (z) , \q n\in{\Bbb Z}, \q w \in L^\infty (\Bbb T; d {\bf m}).
 \label{eq:BH}\end{equation}
     \end{theorem}

  However  results on unbounded Toeplitz operators are very scarce. We can mention only the paper \cite{Hartman} by P.~Hartman and the recent survey \cite{Sarason} by D.~Sarason; see also references in these articles.

 \medskip

{\bf 1.2.}
In this paper, we consider semibounded Toeplitz operators in the space $\ell^2({\Bbb Z}_{+})$.
{\it We always suppose that $t_{n}= \ov{t_{-n}}$ so that the quadratic form \e{eq:QFq} is real and assume that
    \begin{equation}
t[g,g] \geq \gamma \| g\|^2  ,\q g\in{\cal D}, \q \| g\| =\| g\|_{\ell^2({\Bbb Z}_{+})},
 \label{eq:T1}\end{equation}
 for some } $\gamma\in \Bbb R$.
In this   case, we are tempted to define $T$ as a self-adjoint operator corresponding to the quadratic form $t[g,g]$. Such an   operator exists if  the form $t[g,g]$ is closable in the space $\ell^2 ({\Bbb Z}_{+})$, but as is well known this is not always true (see  Example~\ref{ex}, below). We refer to the book \cite{BSbook} for  basic information concerning these notions; they are also briefly discussed in Subsection~2.1.    We recall that, by  definition, the operator corresponding to the form $t[g,g] + \d \| g\|^2 $ is given by the equality $T_{\d} = T+\d I$ (observe that the identity operator $I$ is a Toeplitz operator).
Also by  definition, if a form $t[g,g] $  is closable, then all forms $t[g,g] + \d \| g\|^2$  are closable. Therefore we can suppose that the number $\gamma$ in \e{eq:T1} is positive; for definiteness, we choose $\gamma=1$.
 
  We
proceed from the following well known result (see, e.g., \S 5.1 of the book \cite{AKH}) that is a consequence of the F.~Riesz-Herglotz theorem.

  \begin{theorem}\label{R-H} 
   The condition 
      \begin{equation}
 \sum_{n,m\geq 0} t_{n-m} g_{m}\bar{g}_{n}\geq 0, \q \forall g\in \cal D,
 \label{eq:QF}\end{equation}
  is satisfied if and only if there exists a
  non-negative $($finite$)$ measure $dM(z)$ on the unit circle $\Bbb T$ such that the coefficients $t_{n }$ admit the representations 
  \begin{equation}
  t_{n} =\int_{\Bbb T} z^{-n} dM (z)  , \q n\in{\Bbb Z}.
 \label{eq:WH}\end{equation}
    \end{theorem}
    
    Equations \e{eq:WH} for the measure $dM (z)$ are known as the trigonometric moment problem. 
       Of course  their solution is   unique.  Note that for the Lebesgue measure $d {\bf m} (z)$ we have $t_{0 }=1$ and $t_{n }=0$ for $n\neq 0$. Therefore the measure corresponding to the form
    $t[g,g] + \d \| g\|^2 $ equals $dM  (z)+ \d d {\bf m} (z)$. So we have a one-to-one correspondence between Toeplitz  quadratic forms satisfying estimate  \e{eq:T1} and real measures satisfying the condition $M  (X)\geq \gamma    {\bf m} (X)$ for all Borelian sets $X\subset{\Bbb T}$.
    
   Our goal is to find necessary and sufficient conditions for the form $t[g,g]$ to be  closable. The answer to this question is strikingly  simple.
   
    \begin{theorem}\label{T1} 
    Let the form $t[g,g]$ be given by formula \e{eq:QFq} on elements $g\in {\cal D}$, and let the condition 
 \e{eq:T1} be satisfied. 
   Then the form $t[g,g] $ is closable in the space $\ell^2 ({\Bbb Z}_{+})$ if and only if the measure $dM  (z)$ in   the equations \e{eq:WH} is absolutely continuous.
    \end{theorem}
    
    Of course  Theorem~\ref{T1} means that $dM  (z)=  w(z) d{\bf m}(z)$ where the function $w \in L^1 (\Bbb T; d{\bf m})$ and $w (z)\geq \gamma$. Thus Theorem~\ref{T1} extends Theorem~\ref{B-H} to semibounded operators. The function 
$w (z)$ is known as the symbol of the Toeplitz operator     $T$. So Theorem~\ref{T1}  shows that for a semibounded Toeplitz operator (even defined via the corresponding quadratic form), the symbol exists and is a semibounded function.
   
 \medskip

{\bf 1.3.}   
    The proof of Theorem~\ref{T1} will be given in the next section. The closure of the form $t[g,g] $ is   described in Theorem~\ref{Clo}. 
    
    The discussion of   these results as well as their comparison with  similar statements for Hankel operators are  postponed until Section~3.      We there also explain    shortly how our  results extend to vectorial Toeplitz operators.

  \section{PROOF OF THEOREM~\ref{T1}}  
%

{\bf 2.1.}
Let $t[g,g]$ be a quadratic form defined on a set $\cal D$ dense in a Hilbert space $\cal H$ and satisfying inequality \e{eq:T1} where $\| g\|$ is the norm of  $g\in\cal H$. Suppose that $\gamma=1$,  consider the norm $\| g\|_{T}= \sqrt{t[g,g]}$ and introduce the closure ${\cal D}[t]$ of $\cal D$  in this norm. If ${\cal D}[t]$ can be realized as a subset of $\cal H$, then one says that $t[g,g]$ is closable in the space $\cal H $;  this means that the conditions
      \[
\| g^{(k)}\|\to 0 \q{\rm and}  \q \| g^{(k)} - g^{(l)}\|_{T}\to 0
\]
 as $k,l\to\infty$ imply that $\| g^{(k)}\|_{T}\to 0$. It is easy to see that if $T_{0}$ is a symmetric semibounded operator on $\cal D$, then the form $  t[g,g]=  ( g, T_{0} g)$ is closable.
 
 If the form $t[g,g]$ is closable, then  ${\cal D}[t]\subset \cal H$   is a closed set with respect to the norm $\| \cdot \|_{T}$. In this case $t[g,g]$ is defined by continuity on all $g\in {\cal D}[t]$, and one says that
 the form $t[g,g]$ is closed on   ${\cal D}[t] $.  For a closed form there exists a unique self-adjoint operator $T$ such that $T\geq I$ and
        \begin{align*}
        t[g,h] &=  ( g, T h), \q \forall g\in {\cal D}[t], \q \forall h\in  {\cal D} (T)\subset {\cal D}[t],
        \\
        t[g,g]&= \| \sqrt{T} g\|^2, \q \forall g\in  {\cal D}(\sqrt{T}) = {\cal D}[t].
            \end{align*}
    Note that the domain $ {\cal D} (T) $ of $T$  does not admit an efficient description.

 We are  going to use these general definitions for the space ${\cal H}  = \ell^2 ({\Bbb Z}_{+})$ and the Toeplitz quadratic forms \e{eq:QFq}.
 
 Of course  quadratic forms, in particular, the  Toeplitz  forms, are not necessarily  closable.
 
     \begin{example}\label{ex}
Let $t_{n}=1$ for all $n\in {\Bbb Z}$. Adding the term $\| g\|^2$,  we obtain the form
    \[
t[g,g]= \big|  \sum_{n\geq 0} g_{n}\big|^2 +  \sum_{n\geq 0}| g_{n}|^2 
 \]
satisfying inequality \e{eq:T1}  with $\gamma=1$.
 Define the sequence $g^{(k)}\in {\cal D}$ by the equalities $g_{n}^{(k)}=k^{-1}$ for $0\leq n < k$ and $g_{n}^{(k)}=0$ for $ n\geq k$. Then  $\| g^{(k)}\|=k^{-1/2}\to 0$ as $k\to\infty$. Since  $\sum_{n \geq 0} g_{n}^{(k)}=1$, we  have $ \| g^{(k)} - g^{(l)}\|_{T}= \| g^{(k)} - g^{(l)}\|\to 0$ as $k,l\to\infty$.  Nevertheless  $ \| g^{(k)} \|_{T} \geq 1$.
    \end{example}
    
    Note that the measure $d M(z)$ corresponding to the sequence $t_{n}=1$, $\forall n\in {\Bbb Z}$, is supported by the point $1 \in {\Bbb T}$: $M (\{1\})=1$, $M ({\Bbb T} \setminus \{1\})=0$.
    
    On the other hand, we have  the following simple assertion.
    
       \begin{lemma}\label{ex1}
If a sequence $\{t_{n}\}_{n\in{\Bbb Z}}\in\ell^2 ({\Bbb Z})$, then the form \e{eq:QFq}  is closable.
    \end{lemma}
    
\begin{pf}
   Now we have $  t[g,g]=  ( g, T_{0} g)$ where the symmetric  operator $T_{0}$ is defined by formula \e{eq:HD}  on the set $\cal D$.
    \end{pf}
        
  \medskip
  
  {\bf 2.2.}  
  As  already mentioned, 
  by the proof of Theorem~\ref{T1}  we may suppose that estimate    \e{eq:T1}  is true for $\gamma=1$. According to Theorem~\ref{R-H} the equations  \e{eq:WH} are satisfied with a measure $dM (z)$ such  that $M(X)\geq {\bf m} (X)$ for all Borelian sets  $X\subset \Bbb T$;  in particular, the measure $dM (z)$  is positive. 
  
        Our proof relies on the following auxiliary construction. Let $   L^2 (\Bbb T; dM)$ be the space of functions $u (z)$ on $\Bbb T$ with the norm 
   \[
   \| u\|_{ L^2 (\Bbb T; dM)}=\sqrt{\int_{\Bbb T} | u(z)|^2 dM(z)}.
   \]
          We put
   \begin{equation}
({\cal A}g) (z)=\sum_{n=0}^\infty   g_{n} z^n
 \label{eq:A}\end{equation}
 and      observe that  ${\cal A} g\in  L^2 (\Bbb T; dM)$ for all $g \in\cal D$.
Therefore we can  define   an   operator $A \colon \ell^2 ({\Bbb Z}_{+})\to    L^2 (\Bbb T; dM)$  on domain $\cal D (A)=\cal D$ by the formula $Ag={\cal A}g$.
 In view of  equations  \e{eq:WH}, the form       \e{eq:QFq}   can be written as
   \begin{multline}
   t [g,g] =\sum_{n,m\geq 0}\int_{\Bbb T}z^{-n+m}  d M(z) g_m \ov{g_{n} } 
   \\
     =\int_{\Bbb T}  | ({\cal A} g) (z)|^2 dM(z) =  \| Ag\|^2_{ L^2 (\Bbb T; dM)} , \q g\in \cal D.
 \label{eq:A2}\end{multline}
 This yields the following result.
 
    \begin{lemma}\label{de}
 The form $   t [g,g]$ defined   on $\cal D$   is closable in the space $\ell^2 ({\Bbb Z}_{+})$  if and only if the operator $A \colon \ell^2 ({\Bbb Z}_{+})\to    L^2 (\Bbb T; dM)$ defined   on the  domain $\cal D (A)=\cal D$ 
  is closable.  
    \end{lemma}
     
      Our next goal is to construct the adjoint operator $A^*$.  Observe that       
 for an arbitrary $u\in L^2 (\Bbb T; dM)$, all the integrals
       \begin{equation}
 \int_{\Bbb T}  u(z) z^{-n} dM(z)=: u_{n}, \q n\in {\Bbb Z}_{+},
 \label{eq:A1}\end{equation}
 are absolutely convergent and the sequence $ \{u_{n}\}_{n=0}^\infty$ is bounded. 
 We denote by   ${\cal D}_{*}\subset   L^2 (\Bbb T; dM)$ the set of all  $u \in  L^2 (\Bbb T; dM)$ such that   $\{u_{n}\}_{n=0}^\infty \in  \ell^2 ({\Bbb Z}_{+})$.

   \begin{lemma}\label{LTM}
 The operator $A ^*$ is given by the equality 
      \begin{equation}
        (A^ {*}u)_{n}=  \int_{\Bbb T}  u(z) z^{-n} dM(z), \q n\in {\Bbb Z}_{+},
 \label{eq:a1}\end{equation}
  on the domain ${\cal D}(A^*) ={\cal D}_{*}$. 
    \end{lemma}
    
     \begin{pf}
     Obviously, for all $g\in \cal D$ and all $u\in L^2 (\Bbb T; dM)$, we  have the equality
       \begin{equation}
 (Ag,u)_{L^2 (\Bbb T; dM)} =  \int_{\Bbb T}   \sum_{n=0}^\infty g_{n} z^{n} \ov{u(z)}  dM(z) =   \sum_{n=0}^\infty g_{n}  \bar{u}_{n}
\label{eq:Y1}\end{equation} 
where the sequence $u_{n}$ is defined by relation \e{eq:A1}. If $u\in {\cal D}_{*}$, then
     the right-hand side here equals $(g, A^* u)$. It follows that $   {\cal D}_{*}\subset {\cal D}(A^*)  $.

         Conversely, if $u \in {\cal D}(A^*) $, then 
         \[
         |(Ag, u)_{L^2 (\Bbb T; dM)} | =   |(g, A^* u)_{ \ell^2 ({\Bbb Z}_{+})} | \leq   \| A^* u\|_{ \ell^2 ({\Bbb Z}_{+})}  \,  \| g\|_{ \ell^2 ({\Bbb Z}_{+})}
             \]
          for all $ g \in {\cal D}$. Therefore it follows from equality \e{eq:Y1} that
               \[
\big|   \sum_{n=0}^\infty g_{n}  \bar{u}_{n} \big| \leq  \| A^* u\|_{ \ell^2 ({\Bbb Z}_{+})}  \| g\|_{ \ell^2 ({\Bbb Z}_{+})}, \q \forall g \in {\cal D}.
\]
Since $ {\cal D}$ is dense in $ \ell^2 ({\Bbb Z}_{+})$, we  see that
 $\{ u_{n}\}_{n=0}^\infty \in \ell^2 ({\Bbb Z}_{+})$, and hence $u\in {\cal D}_{*}$.
          Thus   ${\cal D}(A^*)   \subset {\cal D}_{*}$.
          \end{pf}

      Recall that    an operator $A $   is closable if and only if its adjoint operator $$A^* \colon L^2 (\Bbb T; dM) \to \ell^2 ({\Bbb Z}_{+}) $$ is densely defined. We use the notation $ \cl{\cal D}_{*}$ for the closure of the set ${\cal D}_{*}$
      in the space $L^2 ({\Bbb T}; dM)$. So we have obtained an intermediary result.
      
       \begin{lemma}\label{adj}
 The operator $A  $ and the form $t[g,g]$ are closable if and only if  
   \begin{equation}
 \cl{\cal D}_{*} =L^2 ({\Bbb T}; dM). 
 \label{eq:D}\end{equation}
    \end{lemma}

 \medskip
 
 {\bf 2.3.}
Next, we use the   Riesz Brothers theorem. We state it in a  slightly more general form  than in most textbooks.

         \begin{theorem}\label{brothers}
       For a  complex $($finite$)$ measure $d\mu (z)$    on the unit circle $\Bbb T$, put
              \[
 \wh{\mu} (n)=\int_{\Bbb T}  z^{-n}d\mu(z) 
\]
 and suppose that 
 \begin{equation}
\wh{\mu} (n)\in \ell^2 ({\Bbb Z}_{+}). 
 \label{eq:brbr}\end{equation}
     Then  the measure $d\mu   (z)$ is absolutely continuous.
    \end{theorem}
    
    Indeed, in view of \e{eq:brbr}  the function
       \begin{equation}
    f(z):= \sum_{n=0}^\infty  \wh{\mu} (n) z^n
 \label{eq:br1}\end{equation}
 belongs to $L^2 ({\Bbb T}; d{\bf m})$.
    Let us consider an auxiliary measure
    \begin{equation}
    d\mu_{0}(z)= d \mu(z) -  f(z) d {\bf m} (z),
 \label{eq:br2}\end{equation}
 and let
 \[
    \wh{\mu}_{0} (n) =  \wh{\mu} (n) - \int_{\Bbb T} z^{-n}  f(z) d {\bf m} (z),\q n\in {\Bbb Z}_{+}  ,
    \]
    be its Fourier coefficients.
 It follows from \e{eq:br1} that $ \wh{\mu }_{0} (n)= 0$ for all $n\geq 0$. So, by the standard version of the  Riesz Brothers  theorem  (see, e.g., \cite{Hof}, Chapter~4), the measure $  d\mu_{0}(z)$  is absolutely continuous. In view of \e{eq:br2}, the same is true for the measure $  d\mu(z)$.
    
    The following assertion  is almost obvious.
    
      \begin{lemma}\label{ac}
      Suppose that a set ${\cal D}_{*}$ satisfies condition  \e{eq:D}. Let the measures $u(z) dM(z)$ be absolutely continuous
      for all $u\in{\cal D}_{*}$.  Then the measure $  dM(z)$  is also  absolutely continuous.
    \end{lemma}

  \begin{pf}  
  Denote by $\1_{X}$ the characteristic function of a Borelian set $X\subset \Bbb T$.   It follows from \e{eq:D} that there exists a sequence $u_{n}\in{\cal D}_{*}$ such that
  \[
\lim_{n\to\infty}  \| u_{n} -\1_{X} \|_{L^2 ({\Bbb T}; d M)}=0
\]
and hence  
 \[
\lim_{n\to\infty}  \int_{X} u_{n}(z) dM(z)=  \int_{X}   dM(z)= M(X).
\]
If ${\bf m} (X)=0$ and the measures $u_{n}(z) dM(z)$ are absolutely continuous, then the integrals on the left-hand side are zeros so that $M  (X)=0$. 
 \end{pf} 
 
 Now it is easy to conclude the ``only if " part of  Theorem~\ref{T1}. 
  Suppose that the form $t[g,g] $  is closable. Then by Lemma~\ref{adj} the condition  \e{eq:D} is satisfied. 
  By the definition of the set ${\cal D}_{*}$, the Fourier coefficients of  the measures $\mu(z)= u(z) dM(z)$ belong to $\ell^2 ({\Bbb Z}_{+})$. Therefore it follows from   Theorem~\ref{brothers} that   these measures   are absolutely continuous for all $u\in{\cal D}_{*}$.    Hence   by Lemma~\ref{ac} the measure $  dM(z)$  is also absolutely continuous.



  \medskip
  
  {\bf 2.4.}
  It remains to check the converse statement. 
   Let  us consider an auxiliary form 
   \begin{equation}
{\bf t}[g,g] =\sum_{n,m\in {\Bbb Z}} t_{n-m} g_{m}\bar{g}_{n}   
 \label{eq:QFQ}\end{equation}
 in the space $\ell^2({\Bbb Z} )  $ on the set ${\sf D} $ of sequences $g=\{g_{n}\}_{n\in {\Bbb Z}}$  with only a finite number of non-zero components. It is easy to see that the conditions \e{eq:T1} and 
    \begin{equation}
{\bf t}[g,g] \geq \gamma \| g\|^2_{\ell^2({\Bbb Z} )},\q g\in{\sf D},  
 \label{eq:T1x}\end{equation}
 are equivalent. We again suppose that  $\gamma=1$.  Then  the coefficients $t_{n}$  are given by  formula \e{eq:WH} where $M(X)\geq {\bf m} (X)$ for all $X\subset {\Bbb T}$.
  Let  the operator ${\sf A}\colon \ell^2({\Bbb Z}) \to L^2 ({\Bbb T}; dM)$ be defined (cf. \e{eq:A}) by the formula
    \begin{equation}
( {\sf A} g) (z)= \sum_{n\in {\Bbb Z} } z^{n}    g_n=: f(z), \q {\cal D} ({\sf A}  )= {\sf D}.
 \label{eq:A4}\end{equation}
  Quite similarly to \e{eq:A2}, we find that for $g\in \sf D$,
     \begin{equation}
   {\bf t} [g,g] =  \int_{\Bbb T}   | ({\sf A} g) (z)|^2 dM (z) =\| {\sf A} g \|^2_{L^2 ({\Bbb T}; dM)}  . 
 \label{eq:A3}\end{equation}
 
  
  For the completeness of our presentation, let us check that the measure $dM(z)$ is absolutely continuous if the form $   {\bf t} [g,g]$ is closable. This can be done similarly to the proof  in the previous subsection of the same fact for the form $   t [g,g]$, but now  $n\in{\Bbb Z}$ in all formulas and the Riesz Brothers theorem is not needed.   The   operator ${\sf A}^* \colon  L^2 ({\Bbb T}; dM) \to \ell^2 ({\Bbb Z} ) $  acts again by the formula \e{eq:a1}, and it is defined on the set  ${\sf D}_{*} $  of all $u\in  L^2 ({\Bbb T}; dM) $ such that ${\sf A}^* u \in \ell^2 ({\Bbb Z} ) $.
  This means that $u\in  {\sf D}_{*} $ if and only if $u(z) dM(z)= \varphi (z) d {\bf m}(z)$ for some $\varphi \in L^2 ({\Bbb T}; d {\bf m})$. The form $   {\bf t} [g,g]$ is closable if and only  if
   \[
 \cl{\sf D}_{*} =L^2 ({\Bbb T}; dM). 
 \]
Hence Lemma~\ref{ac} implies that the measure $dM(z)$
  is absolutely continuous.

  Next, we show that, for absolutely continuous measures $dM(z) $,  the forms ${\bf t} [g,g]$ are closable.  Now we have  $dM(z)= w(z) d {\bf m} (z)$ where $w \in L^1 ({\Bbb T}; d {\bf m})$. Therefore it follows from \e{eq:A3}   that
   \[
    {\bf t} [g,g ] = {\bf s} [{\sf A}g ,{\sf A} g]  
\]
 where
  \begin{equation}
   {\bf s} [f,f] =  \int_{\Bbb T} w(z) | f (z)|^2 d {\bf m} (z) . 
 \label{eq:A5}\end{equation}
 Since the operator ${\sf A}: \ell^2 ({\Bbb Z})\to L^2 ({\Bbb T};  d {\bf m})$ is unitary, the form ${\bf t} [g,g]$    defined on $\sf D$ is closable in $\ell^2 ({\Bbb Z})$ if and only if the form ${\bf s} [f,f]$   defined on the quasi-polynomials \e{eq:A4}  is closable
in $L^2 ({\Bbb T};  d {\bf m})$. Clearly,   ${\bf s} [f,f]$  is the quadratic form of the operator of multiplication by $ w(z)$.
So it is closable because    $ w\in L^1 ({\Bbb T} ;  d {\bf m})$. Moreover,       ${\bf s} [f,f]$  is closed on the set of all  $f\in L^2({\Bbb T} ;  d {\bf m})$ such the integral \e{eq:A5} is finite. 

Let us summarize the results obtained for the form $  {\bf t} [g,g]$.
 
   \begin{proposition}\label{bro}
Let  the form ${\bf t}[g,g] $ be defined by the equality  
 \e{eq:QFQ}  in the space $\ell^2({\Bbb Z} )  $ on the set ${\sf D} $ of sequences $g= \{g_{n}\}_{n\in {\Bbb Z}}$  with only a finite number of non-zero components.  Let  inequality \e{eq:T1x}  hold, and let the measure $dM(z)$ be defined by 
 the  relation \e{eq:WH}. Then  the form ${\bf t}[g,g] $  is closable   if and only if $dM(z)= w(z) d {\bf m} (z)$ for some $w\in L^1 ({\Bbb T})$, $w(z)\geq \gamma$. In this case
 the closure of ${\bf t}[g,g] $  is given by relations \e{eq:A4}, \e{eq:A3} on the set of all   sequences $ g \in \ell^2({\Bbb Z} )$  such that the integral \e{eq:A3} is finite.
    \end{proposition}

     \medskip

{\bf 2.5.}
  Let us return to Theorem~\ref{T1}.
Obviously, if the form ${\bf t}[g,g]$  is closable   in the space  $\ell^2({\Bbb Z} )  $, then the same is true for the form $t [g,g]$  in the space  $\ell^2({\Bbb Z}_{+} )  $. So if the measure $dM(z)$  is absolutely continuous, then by Proposition~\ref{bro}   both forms ${\bf t}[g,g]$ and $t[g,g]$  are closable. This   concludes the proof of Theorem~\ref{T1}.

\begin{remark}\label{broy}
$\rm(i)$
 For the proof of Theorem~\ref{T1} we used only the ``if " part of  Proposition~\ref{bro}.
On the other hand, the ``only if " part of  Proposition~\ref{bro} is a consequence of  the ``only if " part of Theorem~\ref{T1} proven in Subsection~2.3.   

$\rm(ii)$
  Comparing Theorem~\ref{T1} and Proposition~\ref{bro}, we see that if the  form $t[g,g]$  is  closable, then  the same is true for the form ${\bf t}[g,g]$. As   already noted, the converse assertion is evident.
    \end{remark}

     \medskip
     

{\bf 2.6.}
We now suppose that $dM(z)= w(z) d{\bf m} (z)$ where $w\in L^1 ({\Bbb T}; d {\bf m})$ and $w(z)\geq 1$ so that the form $t[g,g] $  is closable  in the space $\ell^2 ({\Bbb Z}_{+})$.  To describe its closure,  we need a mild additional assumption on $w(z)$. We suppose that the function $w(z)$ is a  Muckenhoupt weight; see, e.g., \S B 5.7 of the book \cite{NK} for various definitions of this notion. One of them is given by the condition 
   \begin{equation}
  \sup_{X \subset {\Bbb T}} {\bf m}(X)^{-2}\int_{X} w (z) d{\bf m} (z) \int_{X} w (z)^{-1} d {\bf m} (z)<\infty
 \label{eq:Mu}\end{equation} 
 where $X$ runs over all subarcs of $\Bbb T$. Let $P_{+}$ be the orthogonal projection of $L^2 ({\Bbb T}; d {\bf m})$
 onto the Hardy class $H^2 ({\Bbb T}; d {\bf m})$ of functions analytic in the unit disc.
 The operator $P_{+}\colon L^2 (\Bbb T; d M) \to L^2 (\Bbb T; dM)$ is bounded if and only if $w (z)$ is a Muckenhoupt weight.  
 Recall that   the operator  ${\cal A} $ is defined by formula  \e{eq:A}. Obviously, ${\cal  A}g\in H^2 ({\Bbb T};  d {\bf m}) $ for all $g\in \ell^2 ({\Bbb Z}_{+})$.

\begin{theorem}\label{Clo} 
     Let    the coefficients $t_{n}$ be given by formula \e{eq:WH} where $dM(z)= w(z) d {\bf m}(z)$,  $w\in L^1 ({\Bbb T}; d {\bf m})$ and $w (z)\geq 1$. Suppose that $w(z)$ is a Muckenhoupt weight.  
       Then  the closure of the form $t[g,g]$ defined on  ${\cal D}$ by  \e{eq:QFq}
 is given by the equality 
  \begin{equation}
   t [g,g]     =\int_{\Bbb T}  | ({\cal A} g) (z)|^2 dM(z)  
 \label{eq:A2x}\end{equation}
  on the set ${\cal D}[t]$ of all $g\in \ell^2 ({\Bbb Z}_{+})$ such that the right-hand side of \e{eq:A2x} is finite.
        \end{theorem}
        
 \begin{pf}
        Observe that the operator ${\cal A} : \ell^2 ({\Bbb Z}_{+})\to H^2 (\Bbb T; d {\bf m})$  is unitary and  ${\cal A}  {\cal D} =:  {\cal P}$ consists of all polynomials \e{eq:A}. Let $\cl {\cal P}$ be   the closure of $ {\cal P}$ in $ L^2 (\Bbb T; d M)$. So the assertion of Theorem~\ref{Clo} is equivalent to the equality
          \begin{equation}
\cl {\cal P}= H^2 (\Bbb T;  d {\bf m}) \cap L^2 (\Bbb T; d M). 
 \label{eq:A2x1}\end{equation}
 Since the convergence in $ L^2 (\Bbb T; d M)$ is stronger than that in $ L^2 (\Bbb T; d{\bf m})$, we have $\cl {\cal P} \subset H^2 (\Bbb T; d {\bf m})$ and therefore  the left-hand side of \e{eq:A2x1} is contained in its right-hand side.
 
 It remains to prove the opposite inclusion. Recall that if all Fourier coefficients of some complex measure on $\Bbb T$ are zeros, then this measure is also zero. Suppose that
   $u\in L^2 (\Bbb T; d M)$ is orthogonal in $ L^2 (\Bbb T; d M)$ to the functions $z^n$  for all $n\in{\Bbb Z}$. Then applying the above fact to the measure $u dM$, we see that $u=0$. Therefore
    quasi-polynomials $f (z)=\sum_{n\in{\Bbb Z}} a_{n  } z^n$ (the sum consists of a finite number of terms) are dense in $   L^2 (\Bbb T;  d M)$. So for every  $u\in   L^2 (\Bbb T;  d M)$ there exists  a sequence of quasi-polynomials $f_{k} (z)$ such that
      \[
\lim_{k\to\infty} \| u- f_{k} \|_{L^2 (\Bbb T; d  M)}=0, \q dM=w d{\bf m}. 
 \]
 Since  $w(z)$ is a Muckenhoupt weight, this implies that 
   \[
\lim_{k\to\infty} \| P_{+}u- P_{+} f_{k} \|_{L^2 (\Bbb T; dM)}=0. 
\]
 So if $u=P_{+}u$ and $\varphi_{k} =P_{+} f_{k} \in {\cal P}$, then $\varphi_{k}  \to u$ as $k\to\infty$  in $L^2 (\Bbb T; d M)$ which also implies the convergence in $H^2 (\Bbb T;  d{\bf m})$.
  \end{pf}
  
  \begin{remark}\label{Clo1} 
     Recall that the operator $A$ was defined by formula \e{eq:A}  on the domain ${\cal D} (A)= {\cal D} $.
     Of course its closure $\cl{A}=A^{**}$.  Let $A_{\rm max}$ be given by the formula $A_{\rm max}g= {\cal A}g$   on the domain $\cal D(A_{\rm max})$ that consists of all $g\in \ell^2 ({\Bbb Z}_{+})$ such that ${\cal A}g\in L^2 (\Bbb T; dM)$. Then  the assertion of Theorem~\ref{Clo}  is equivalent to the equality
     \[
     \cl{A}= A_{\rm max}.
     \]
        \end{remark}

\section{DISCUSSION}  

{\bf 3.1.}
Let us state a consequence of Theorem~\ref{T1} in terms of the entries $t_{n}$ of the form \e{eq:QFq}.

  \begin{proposition}\label{T1m} 
Suppose that the condition 
 \e{eq:T1} is satisfied. 
If  the form $t[g,g] $  is closable  in the space $\ell^2 ({\Bbb Z}_{+})$, then  $t_{n}\to 0$  as $|n|\to\infty$. If $ \{t_n\}_{n\in {\Bbb Z}}\in \ell^2 ({\Bbb Z})$, then the form $t[g,g] $  is closable in the space $\ell^2 ({\Bbb Z}_{+})$.
    \end{proposition}
    
     \begin{pf}
   If  the form $t[g,g] $  is closable, then, by   Theorem~\ref{T1}, the measure $dM(z)$ in the representation
   \e{eq:WH}  is absolutely continuous. Therefore its Fourier coefficients  $t_{n}\to 0$  as $| n| \to\infty$.
   Conversely,  if $  \{t_n\}_{n\in {\Bbb Z}} \in \ell^2 ({\Bbb Z})$, then $d M(z)= w (z) d {\bf m} (z)$ with $w \in L^2 ({\Bbb T}; d {\bf m} )$. Therefore, again  by   Theorem~\ref{T1},  the form $t[g,g] $  is closable  (this result was already stated in Lemma~\ref{ex1}).
        \end{pf}
        
        There is a  gap between necessary  and sufficient conditions on $ t_n$ in Proposition~\ref{T1m}. Apparently it cannot be significantly reduced. Recall that  by the Wiener theorem (see, e.g., Theorem~XI.114 in \cite{RS}), if the Fourier coefficients $ t_n$ of some measure $dM(z)$  tend to zero, then this measure is necessarily continuous, but it may be singular with respect to the Lebesgue measure.   Thus  the condition        $ t_n\to 0$ as $|n |\to \infty$ does not imply that the   measure $dM(z)$ defined by \e{eq:WH} is absolutely continuous. So in accordance with  Theorem~\ref{T1}    
   the   corresponding Toeplitz quadratic form $t[g,g]$ may be unclosable.
        
          Astonishingly, the sufficient condition  $ \{t_n\}_{n\in {\Bbb Z}}\in \ell^2 ({\Bbb Z})$  for  the absolute continuity of the measure $dM(z)$ turns out to be very sharp.  Indeed,  for every $l\in{\Bbb Z}_{+}$, O.~S.~Iva\v{s}\"{e}v-Musatov constructed in \cite{I-M}  a singular measure such that its Fourier coefficients satisfy the estimate
  \[
  t_{n} =O \big( ( n \ln n \ln\ln n\cdots \ln_{(l)} n)^{-1/2}  \big)
  \]
  (here $ \ln_{(l)} n$ means that the logarithm  is applied $l$ times to $n$).
  Examples of singular  continuous measures of such type go back to D.~E.~Menchoff \cite{Men}. A comprehensive survey of various constructions of singular  continuous measures  with decaying Fourier coefficients can be found in  \cite{Brow}.

\medskip

{\bf 3.2.}
There is a certain parallelism between the theories of Toeplitz and Hankel operators. For example, the criteria of  boundedness
of Toeplitz and of Hankel operators due to Toeplitz (see Theorem~\ref{B-H}) and to Nehari \cite{Nehari}, respectively, look formally similar. Toeplitz   quadratic forms are linked to the trigonometric  moment problem while Hankel quadratic forms are linked to the  power moment problem. The following result obtained by Hamburger  in \cite{Hamb} plays the role of Theorem~\ref{R-H}.

   \begin{theorem}\label{Hamb} 
   The condition 
      \begin{equation}
 q[g,g]=\sum_{n,m\geq 0} q_{n+m} g_{m}\bar{g}_{n}\geq 0, \q \forall g\in \cal D,
 \label{eq:QFh}\end{equation}
  is satisfied if and only if there exists a
  non-negative measure $d{\sf M}(x)$ on $\Bbb R$ such that the coefficients $q_{n }$ admit the representations 
    \begin{equation}
q_{n } = \int_{-\infty}^\infty x^n d {\sf M}(x), \q \forall n=0,1,\ldots.
 \label{eq:WHh}\end{equation}
    \end{theorem}

It is interesting to compare  Theorem~\ref{T1}   with the corresponding result for
 Hankel operators.    Let us state necessary  and sufficient conditions guaranteeing that a   Hankel quadratic form is closable. 
 
   \begin{theorem}\label{Hamb1}\cite[Theorem 1.2]{Ya}   
     Let   assumption \e{eq:QFh} be satisfied.  Then the following conditions are equivalent: 
     
     \begin{enumerate}[\rm(i)]
\item
The form $q[g,g]$ defined on $\cal D$  is closable   in the space $\ell^2 ({\Bbb Z}_{+})$.

 \item
The matrix elements  $ q_n\to 0$ as $n\to \infty$.

\item
The measure $d {\sf M}  (x)$ defined by equations  \e{eq:WHh} satisfies  the condition
  \[
{\sf M} ({\Bbb R}\setminus (-1,1) )=0
\]
$($to put it differently, $\supp {\sf M}\subset [-1,1]$ and ${\sf M}(\{-1\}) = {\sf M}(\{1\})=0)$.
\end{enumerate}
        \end{theorem}
        
  Thus in contrast to Toeplitz  quadratic forms, the condition        $ q_n\to 0$ as $n\to \infty$ is necessary  and sufficient for a   Hankel quadratic form \e{eq:QFh} to be   closable. 
  
  For Hankel quadratic forms, an analogue of Theorem~\ref{Clo} (see Theorem~3.4 in \cite{Ya}) is true without additional assumptions on the measure
  $d {\sf M}  (x)$,  but its proof requires  substantial work.

\medskip

{\bf 3.3.}
Theorem~\ref{T1}  automatically extends to vectorial Toeplitz operators. In this case $g=\{g_{n}\}_{n\in{\Bbb Z}_{+}}$ where  $g_{n}$ are elements of some auxiliary Hilbert space $\mathfrak{N}$, and $t_{n}$ are bounded operators in $\mathfrak{N}$. We now suppose that $t_{n}= t_{-n}^*$ and
 \begin{equation}
t[g,g] =\sum_{n,m\geq 0} \la t_{n-m} g_{m}, g_{n}   \ra_{\mathfrak{N}}  \geq \gamma \| g\|^2_{L^2 ({\Bbb Z}_{+} ;\mathfrak{N})},\q \forall g\in \cal D.
 \label{eq:QFv}\end{equation}
 The vectorial version of Theorem~\ref{R-H}  means that inequality \e{eq:QFv} for $\gamma=0$ is equivalent to the representation \e{eq:WH}   with a non-negative operator valued measure $dM(z)$. Let us state a generalization   of Theorem~\ref{T1} to  the vector case.
 
     \begin{theorem}\label{Tv} 
    Let the condition \e{eq:QFv} be satisfied. 
   Then the form $t[g,g] $ is closable in the space $\ell^2 ({\Bbb Z}_{+};\mathfrak{N})$ if and only if  
     \[
  t_{n} =\int_{\Bbb T} z^{-n} w(z) d{\bf m} (z) , \q n\in{\Bbb Z}, \q w(z)\colon \mathfrak{N}\to\mathfrak{N},
 \]
 where $w(z)\geq \gamma I_{\mathfrak{N}}$ and the operator valued function $w(z)$ belongs to $  L^1 ({\Bbb T}, d{\bf m}  )$.
    \end{theorem}  
    
    Theorem~\ref{Clo}     and its proof  also extend  to the vectorial case provided the projector $P_{+}$ is a bounded operator   in the space  $  L^2 ({\Bbb T}, w d{\bf m}  ;\mathfrak{N})$.   Note that there is a necessary and sufficient condition  (see the paper \cite{T-V}) for the boundedness of this operator generalizing the scalar condition
          \e{eq:Mu}; the result of \cite{T-V} requires however that $\dim \mathfrak{N}<\infty$.
          
          \medskip
          
          The author thanks G.~Rozenblum for a useful discussion. 
 

\end{document}